# ON A GENERALIZATION OF OBRESHKOFF–EHRLICH METHOD FOR SIMULTANEOUS EXTRACTION OF ALL ROOTS OF POLYNOMIALS OVER AN ARBITRARY CHEBYSHEV SYSTEM

## A. Iliev, Khr. Semerdzhiev

Let functions $j_0(x),\ldots,j_n(x)$ be smooth enough and constitute a Chebyshev system on the interval $(a,b)$. Let us consider a generalised polynomial on that system:

$$(1)\quad f(x) = \sum_{j=0}^{n} a_j j_j(x),$$

whose zeros $x_1,\ldots,x_m$ are of multiplicities $a_1,\ldots,a_m$ respectively. Assume $\sum_{j=1}^{m} a_j = n$.

Denote the approximations of the $k^{\text{th}}$ iteration zeros of (1) with $x_1^{[k]},\ldots,x_m^{[k]}$.

The polynomial $f(x)$ can be represented up to a multiplier constant in the following determinant form:

$$(2)\quad f(x) = \begin{vmatrix} j_0(x) & j_1(x) & \ldots & j_n(x) \\ j_0(x_1) & j_1(x_1) & \ldots & j_n(x_1) \\ \ldots & \ldots & \ldots & \ldots \\ j_0^{(a_1-1)}(x_1) & j_1^{(a_1-1)}(x_1) & \ldots & j_n^{(a_1-1)}(x_1) \\ \ldots & \ldots & \ldots & \ldots \\ j_0(x_m) & j_1(x_m) & \ldots & j_n(x_m) \\ \ldots & \ldots & \ldots & \ldots \\ j_0^{(a_m-1)}(x_m) & j_1^{(a_m-1)}(x_m) & \ldots & j_n^{(a_m-1)}(x_m) \end{vmatrix} \overset{\text{def}}{=} \det \begin{bmatrix} j_0 & j_1 & \ldots & j_n \\ x & x_1 & \ldots & x_m \\ a_0 & a_1 & \ldots & a_m \end{bmatrix},\ a_0 \overset{\text{def}}{=} 1.$$

For simultaneous extraction of the roots of equation $f(x)=0$ we will apply the following iteration formula:

$$(3)\quad x_i^{[k+1]} = x_i^{[k]} - a_i f(x_i^{[k]})\left[f'(x_i^{[k]}) - f(x_i^{[k]})Q_i'(x_i^{[k]})\left[(a_i+1)Q_i(x_i^{[k]})\right]^{-1}\right]^{-1},$$
$$i = \overline{1,m},\quad k = 0,1,2,\ldots,$$

where $Q_i(x) = \dfrac{d^{a_i}}{dx^{a_i}} \det \begin{bmatrix} j_0 & j_1 & \ldots & j_n \\ x & x_1^{[k]} & \ldots & x_m^{[k]} \\ a_0 & a_1 & \ldots & a_m \end{bmatrix},\ i = \overline{1,m}$.

Comparing this iteration formula with (13) which we have obtained in our earlier work [3,7], it can be seen that the new formula (3) requires the evaluation of only $f'(x)$ instead of $f^{(a_i)}(x)$ as in (13). Method (3) is simpler than (13) and nevertheless retains its cubic convergence. Besides, the proof of the theorem in [3] is quite complicated. Here we apply a new approach in the theoretical treatment of iterations, which leads to a considerable simplification of the proof of convergence. This formula is a generalization of the iteration formulas of Ehrlich [1] and Obreshkoff [4].

**Definition.** For the functions $p(x_j)$ and $q(x_j)$ we will write $p(x_j) \equiv q(x_j) \mod (x_j - x_j^{[k]})^g, j = \overline{1,m}$ ($p(x_j)$ comparable with $q(x_j)$ modulo $(x_j - x_j^{[k]})^g, j = \overline{1,m}$), when the difference $p(x_j) - q(x_j)$ can be represented as $h(x_j - x_j^{[k]}), j = \overline{1,m}$, where the function $h(y_j)$, ($y_j = x_j - x_j^{[k]}$) whose derivatives of order up to $g-1$ are equal to zero when $y_j = 0$.

At $g = 1$ by expanding the function $h(y_j)$ in a Maclaurin's series to terms of the first order with respect to $y_j$ we obtain

$$h(x_j - x_j^{[k]}) = \sum_{s=1}^{m} \left.\frac{\partial h}{\partial y_s}\right|_{(q\, y_s)} y_s = \sum_{p=1}^{m} A_p (x_p - x_p^{[k]}), j = \overline{1,m}, \quad \text{where the derivatives}$$

$A_p = \frac{\partial h}{\partial y_p}, p = \overline{1,m}$ are calculated at point $(q(x_1 - x_1^{[k]}), \ldots, q(x_m - x_m^{[k]})), 0 < q < 1$.

As an example, we will show that

(4)
$$Q_i(x) \equiv f^{(a_i)}(x) \mod (x_j^{[k]} - x_j), j = \overline{1,m}$$
$$Q_i'(x) \equiv f^{(a_i+1)}(x) \mod (x_j^{[k]} - x_j), j = \overline{1,m}.$$

The second row of $Q_i(x)$ we represent in the following way

$$(j_0(x_1) + [j_0(x_1^{[k]}) - j_0(x_1)], j_1(x_1) + [j_1(x_1^{[k]}) - j_1(x_1)], \ldots, j_n(x_1) + [j_n(x_1^{[k]}) - j_n(x_1)]) =$$
$$(j_0(x_1), j_1(x_1), \ldots, j_n(x_1)) + (j_0'(\mathbf{x}_{0,1}^{0,k})(x_1^{[k]} - x_1), j_1'(\mathbf{x}_{1,1}^{0,k})(x_1^{[k]} - x_1), \ldots, j_n'(\mathbf{x}_{n,1}^{0,k})(x_1^{[k]} - x_1)).$$

From which we see that $Q_i(x) \equiv Q_{i,2}(x) \mod (x_1^{[k]} - x_1)$, $i = \overline{1,m}$. The quantity $Q_{i,2}(x)$ is the determinant $Q_i(x)$ in which the second row is replaced by $(j_0(x_1), \ldots, j_n(x_1))$. Analogously we find that $Q_i(x) \equiv Q_{i,3}(x) \mod (x_1^{[k]} - x_1)$, $i = \overline{1,m}$, where $Q_{i,3}(x)$ is $Q_{i,2}(x)$ in which the third row is replaced by $(j_0'(x_1), \ldots, j_n'(x_1))$. Following this procedure we find that $Q_i(x) \equiv f^{(a_i)}(x) \mod (x_j^{[k]} - x_j), j = \overline{1,m}, i = \overline{1,m}$. Analogously we verify the formula $Q_i'(x) \equiv f^{(a_i+1)}(x) \mod (x_j^{[k]} - x_j), j = \overline{1,m}, i = \overline{1,m}$.

We transform formula (3) to get

(5) $x_i^{[k+1]} - x_i = \Phi(x_i^{[k]})[\Psi(x_i^{[k]})]^{-1}, \quad i = \overline{1,m}$,

where

(6)
$$\Phi(x) = (a_i + 1)[(x - x_i)f'(x) - a_i f(x)]Q_i(x) - (x - x_i)f(x)Q_i'(x)$$
$$\Psi(x) = (a_i + 1)f'(x)Q_i(x) - f(x)Q_i'(x).$$

From (6) it is immediately verifiable that

(7) $\Phi'(x_i) = \Phi''(x_i) = \ldots = \Phi^{(a_i-1)}(x_i) = \Phi^{(a_i)}(x_i) = 0$.

It follows from (4)

(8) $\Phi^{(a_i+1)}(x_i) \equiv 0 \mod (x_j^{[k]} - x_j), j = \overline{1,m}, i = \overline{1,m}$.

Therefore from (8), we have

(9) $\Phi^{(a_i+1)}(x_i) = \sum_{s=1}^{m} A_{is}(x_s^{[k]} - x_s)$, $i = \overline{1,m}$.

It follows from (7) and from Taylor's formula that

(10) $\Phi(x_i^{[k]}) = \Phi^{(a_i+1)}(x_i)\frac{1}{(a_i+1)!}(x_i^{[k]} - x_i)^{a_i+1} + \frac{1}{(a_i+2)!}\Phi^{(a_i+2)}(z_i^{[k]})(x_i^{[k]} - x_i)^{a_i+2}$.

Analogously, we find $\Psi'(x_i) = \ldots = \Psi^{(a_i-2)}(x_i) = 0$ and

(11) $\Psi(x_i^{[k]}) = \frac{1}{(a_i-1)!}\Psi^{(a_i-1)}(x_i)(x_i^{[k]} - x_i)^{a_i-1} + \frac{1}{a_i!}\Psi^{(a_i)}(h_i^{[k]})(x_i^{[k]} - x_i)^{a_i}$, $i = \overline{1,m}$.

We substitute (10) and (11) in formula (5) and obtain

(12) $x_i^{[k+1]} - x_i = \dfrac{\frac{1}{(a_i+1)!}(x_i^{[k]} - x_i)^2 \sum_{s=1}^{m} A_{is}(x_s^{[k]} - x_s) + \frac{1}{(a_i+2)!}\Phi^{(a_i+2)}(z_i^{[k]})(x_i^{[k]} - x_i)^3}{\frac{1}{(a_i-1)!}\left[(a_i+1)[f^{(a_i)}(x_i)]^2 + \sum_{s=1}^{m} B_{is}(x_s^{[k]} - x_s)\right] + \frac{1}{a_i!}\Psi^{(a_i)}(h_i^{[k]})(x_i^{[k]} - x_i)}$,

$i = \overline{1,m}$.

**Theorem.** Let $A \stackrel{\text{def}}{=} \max_{s,i=1,m}\{\sup|A_{is}|\}$, $B \stackrel{\text{def}}{=} \max_{s,i=1,m}\{\sup|B_{is}|\}$, $M \stackrel{\text{def}}{=} \max_{i=1,m}\{\sup|\Phi^{(a_i+2)}(x)|\}$, $K \stackrel{\text{def}}{=} \max_{i=1,m}\{\sup|\Psi^{(a_i)}(x)|\}$, $P \stackrel{\text{def}}{=} \min_{i=1,m}\{\inf|f^{(a_i)}(x_i)|\}$. Let $c > 0$ and $0 < q < 1$ be real numbers such that the following inequalities are fulfilled

$c^2\left[mA + M(a_i+2)^{-1}\right] < (a_i+1)^2 a_i P^2 - c(a_i+1)[a_i Bm + K]$, $i = \overline{1,m}$. Let us assume that $|x_i^{[k]} - x_i| < cq^{3^k}$ for some $k \in N \cup \{0\}$ and $i = \overline{1,m}$. Then $|x_i^{[k+1]} - x_i| < cq^{3^{k+1}}$, $i = \overline{1,m}$.

Thus method (3) has a third order of convergence.

In the particular case when $j_s(x) = x^s$, $s = \overline{1,n}$, i.e. $f(x)$ is of the kind $f(x) = a_n x^n + a_{n-1} x^{n-1} + \ldots + a_0$ we obtain $\dfrac{Q_i'(x_i^{[k]})}{(a_i+1)Q_i(x_i^{[k]})} = \sum_{\substack{j=1\\j\neq i}}^{m} a_j(x_i^{[k]} - x_j^{[k]})^{-1}$, which coincides with the finding in [2,6,7].

**Example.** We will compare the newly found generalization (3) with that obtained by formula [3,7]

(13) $x_i^{[k+1]} = x_i^{[k]} - f^{(a_i-1)}(x_i^{[k]})\left[f^{(a_i)}(x_i^{[k]}) - f^{(a_i-1)}(x_i^{[k]})Q_i'(x_i^{[k]})[2Q_i(x_i^{[k]})]^{-1}\right]^{-1}$,

$i = \overline{1,m}$, $k = 0,1,2,\ldots$

From the base system of functions $\{1, x^2, \sin 3x, \exp(-x), (1+x^2)^{-1}\}$ a generalised polynomial is constructed having $x_1 = -0.5$ and $x_2 = 3$ as its double zeros.

| k | method (3) | | method (13) | |
|---|---|---|---|---|
| | $x_1^{[k]}$ | $x_2^{[k]}$ | $x_1^{[k]}$ | $x_2^{[k]}$ |
| 0 | -0.4000000000 | 2.8000000000 | -0.4000000000 | 2.8000000000 |
| 1 | -0.5001904855 | 2.9812593584 | -0.5021054 | 2.9677106 |
| 2 | -0.5000000001 | 2.9999296686 | -0.500000081 | 2.99935 |
| 3 | -0.5000000000 | 3.0000000000 | -0.5000000000 | 2.9999999915 |

Numerical experiments were performed applying method (3) and showed absolutely identical results (to the $11^{th}$ digit) as the generalization of the Obreshkoff-Ehrlich method for extraction of multiple roots in the particular cases of algebraic, trigonometric, and exponential polynomials [2,7].

Finally we will point out that the method (3) and other methods from this group [5,7,8] for simultaneously extraction of all roots of polynomial equations are quite appropriate for applying to computers with parallel processors.

University of Plovdiv
http://www.pu.acad.bg
Faculty of Mathematics and Informatics
http://www.fmi.pu.acad.bg
24 Tzar Assen Str.
Plovdiv 4000
e-mail: aii@pu.acad.bg
URL: http://anton.iliev.tripod.com